\documentclass{article}
\usepackage{amssymb,amsmath,graphicx,multirow}
 
\newtheorem{theorem}{Theorem}
\newtheorem{corollary}{Corollary}

\newtheorem{lemma}{Lemma}

\title{Partitions with Prescribed Hooksets}

\author{William J. Keith, University of Lisbon\footnote{The first author began work on this paper at Drexel University and concluded while working at the University of Lisbon.}, and Rishi Nath, CUNY-York \\ Author 1 email: wjk150@cii.fc.ul.pt; Author 2 email: rnath@york.cuny.edu}

\begin{document}

\maketitle

\abstract{The combinatorial properties of partitions with various restrictions on their hooksets are explored.  A connection with numerical semigroups extends current results on simultaneous $s/t$-cores.  Conditions that suffice for a partition to possess required hooks are developed.}

\section{Introduction}

Questions concerning the number of partitions with hooksets of various types are of interest in several fields.  Partitions with hook lengths avoiding multiples of $t$, called $t$-cores, are useful in representation theory (where they characterize the $t$-modularly irreducible representations of $S_n$, if $t$ is prime, see \cite{RepSource}) and in standard partition theory, such as for the properties of the crank function \cite{CrankSource}.  The hookset of a partition encodes information about other combinatorial objects associated to that partition, perhaps most famously through the \emph{tableaux number} $\frac{n!}{\prod h_{ij}}$ of a partition $\lambda$, which is the degree of the corresponding irreducible representation of $S_n$ and also counts the number of Young tableaux that have the shape of that partition.

One natural generalization of the $t$-core requires avoiding hooks that are multiples of more than one length; \cite{STCoresPrime} and \cite{STCoresNotPrime} study partitions simultaneously $s$-core and $t$-core, for two numbers $s$ and $t$ coprime and not coprime respectively.  In this article one of the primary results is

\begin{theorem}\label{SCores} Let $S = \{ s_1, \dots , s_k \}$ be any set of positive integers.  The set of partitions which are simultaneously $s_i$-core for all $s_i \in S$ is finite if and only if $gcd(S) = 1$.
\end{theorem}

A simple upper bound for the largest $n$ with such a partition can also be obtained.  These are both consequences of an observation made as Theorem 2.2 in \cite{AddArt}, which can be rephrased as:

\begin{theorem}\label{HooksetMonoids} Among all partitions, the hooksets of length $g$ are exactly the complements of numerical monoids of genus $g$.
\end{theorem}

We give an explicit bijection that proves the theorem and derive the above and a few other consequences in the theory of partitions; it also seems reasonable that this connection could be of use in the theory of numerical monoids.  For example, "maximum possible" hooksets, which are initial segments of the integers, correspond to the previously-studied sequence of numerical sets with "minimum possible" atom monoids:

\begin{corollary}\label{EmptyAtom} The number of partitions with hookset $Hk(\lambda) = \{1,2,\dots,n\}$ equals the number of numerical sets $S$ with atom monoid $A(S) = \{0,n+1,n+2,n+3,n+4,\dots\}$.
\end{corollary}

This immediately raises complementary questions on the combinatorics of hooksets which are required to contain certain values; for example, to be an initial segment of a given sequence $S$.  The best answer might be a counting function once $S$ was described.  This does not seem to be an easy question, though we can tackle it for the extreme of partitions which have very dense hooksets.  A simple abacus argument does establish the sufficient condition

\begin{theorem}\label{DyckHooks} Given an arbitrary set of required hooks $\{s_i\}$, a partition with a profile which lies above the diagonal at points of distance $s_i$ from the upper and lower corners along its profile will contain the hooks $s_i$ in its hookset.
\end{theorem}

So some, though by no means all, of the partitions referred to in Corollary \ref{EmptyAtom} are described thus:

\begin{corollary}\label{StrictCorollary} Strict Dyck paths (partitions with largest part equal to the number of parts, in which the $i$-th smallest part is strictly greater than $i$) contain all hooks up to their largest hook.
\end{corollary}

Regarding the map used to prove Theorem \ref{HooksetMonoids}, it should be noted that a search through the literature does find a very similar bijection used in one paper, \cite{BAdMier}, by Maria Bras-Amor\'{o}s and Anna de Mier.  With some alterations, that map has the same hook property, although this property is not remarked upon in that paper; that paper was interested in the symmetry of the path and the fact that image of a numerical semigroup becomes the profile of a Dyck path.  A more explicit comparison is given after the map is completely described.

\subsection{Notation and Definitions}

The reader familiar with the terminology and notation of partitions, abacus techniques, and numerical sets and semigroups may skip this section.

Let $n$ be a nonnegative integer.  We say that a finite nonincreasing sequence $\lambda$ \emph{partitions} $n$, written $$\lambda = (\lambda_1,\lambda_2,\dots,\lambda_k) \vdash n$$

\noindent if $\lambda_i \in \mathbb{N}$, $\lambda_i \geq \lambda_{i+1} \geq 1$, $\sum \lambda_i = n$.

We draw a partition with its \emph{Young diagram}, which is an array of squares occupying the fourth quadrant with lower right corner at points $(i,-j)$ if $\lambda_j \geq i$.

\begin{center}\begin{tabular}{c} $\Box \Box \Box \Box \Box \Box \Box$\vspace{-0.05in} \\ $\Box \Box \Box \Box \Box \Box$\phantom{$\Box$}\vspace{-0.05in}\\$\Box \Box \Box \Box \Box \Box$\phantom{$\Box$}\vspace{-0.05in}\\$\Box \Box \Box \Box$\phantom{$\Box \Box \Box$}\vspace{-0.05in}\\$\Box \Box \Box$\phantom{$\Box \Box \Box \Box$}\vspace{-0.05in}\\$\Box \Box$\phantom{$\Box \Box \Box \Box \Box$}\vspace{-0.05in}\\$\Box \Box$\phantom{$\Box \Box \Box \Box \Box$}\vspace{-0.05in}\\$\Box$\phantom{$\Box \Box \Box \Box \Box \Box$}\vspace{-0.05in}\\$\Box$\phantom{$\Box \Box \Box \Box \Box \Box$}\end{tabular}\end{center}

The above is the partition $(7,6,6,4,3,2,2,1,1) \vdash 32$.

The \emph{hook} $h_{ij}$ of a partition is the sum of the lengths of the "arm" and "leg" depending from the square representing $(i,-j)$ in the diagram; symbolically, $h_{ij} = \lambda_i - j + |\{\lambda_a \geq j \} |$.

In the example above, the hook $h_{11} = 15$; we will refer to $h_{11}$ as the \emph{outer hook} of a partition.  At the farthest right corner, the hook $h_{71} = 1$.

The \emph{hookset} of a partition will be $Hk(\lambda)=\{h_{ij}\}$; here, the hookset is $Hk(\lambda) = \{1,2,3,4,5,6,7,9,10,12,13,15\}$.  It is to be distinguished from the hook multiset; for example, the partition above has 6 hooks of size 1.

Partitions which have no hooks $h_{ij} = m k$ for any $m$ are called $k$-\emph{cores}.  Some known facts on $k$-cores (see \cite{STCoresNotPrime}, \cite{STCoresPrime}) that will be of use to us are:

\begin{lemma}\label{CoreFacts}
\begin{itemize}
\item If $\lambda$ is a $k$-core, it is a $pk$-core for any $p$, i.e. if $H(\lambda)$ contains $c$ it contains all divisors of $c$.
\item The only 2-cores are given by the staircase partitions $(k,k-1,k-2,...,3,2,1)$.  Hence the only numbers with any partition whose hooks are all odd are $\left( {k \atop 2} \right)$, and the hooksets of these partitions are exactly $\{1,3,5,...,2k-1\}$.
\item Partitions which are both an $s$-core and a $t$-core with $s$ and $t$ coprime form a finite set; the largest member partitions $(s^2-1)(t^2-1)/24$.
\end{itemize}\end{lemma}

Proofs of $k$-core theorems often involve the \emph{abacus} of a partition.  The \emph{profile} of a partition is the set of southmost and eastmost edges of the boxes in its Young diagram:

\begin{center}\begin{tabular}{c}
\includegraphics[scale=0.8]{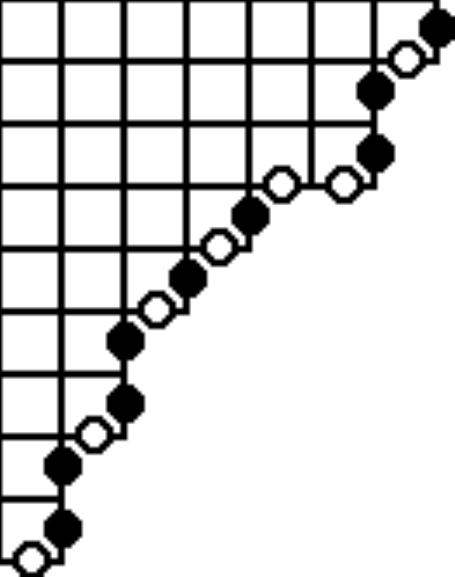}
\end{tabular}\end{center}

South edges, or right steps, are marked with \emph{spacers} while east edges, or up steps, are marked with \emph{beads}.  The list of these forms the \emph{abacus}.

\begin{center}\begin{tabular}{c}
\includegraphics[scale=0.8]{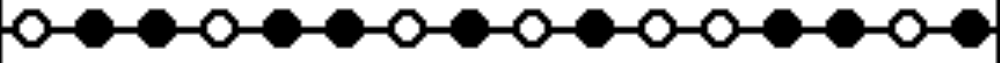}
\end{tabular}\end{center}

The abacus can be wrapped on itself, reading left to right, then top to bottom, to create a useful visual display of the partition's hookset:

\begin{center}
\begin{tabular}{c}
$\circ \bullet \bullet$\vspace{-0.07in}\\
\hline $\circ \bullet \bullet$\vspace{-0.085in}\\
\hline $\circ \bullet \circ$\vspace{-0.085in}\\
\hline $\bullet \circ \circ$\vspace{-0.075in}\\ 
\hline $\bullet \bullet \circ$\vspace{-0.075in}\\
\hline $\bullet \circ \circ$\vspace{-0.075in}\\ \hline
\end{tabular}
\end{center}

The profile starts at the first (top, left) spacer and ends at the last bead (bottom, rightmost); hooks lengths are the distances from any spacer to any beads lower or to the right.  Thus, when the wrapping uses horizontal rods of length $k$, any spacer with a bead directly below it corresponds to a hook of length $k$.  

The abacus of a $k$-core, then, if wrapped on rods of length $k$, must consist of columns which are either empty (i.e. all spacers), or columns of beads that start at the first row and descend without interruption until that column stops, and then are spacers until the profile ends.  Furthermore, the first column must be empty.  We enumerate the columns and rods starting from 0, so that the first spacer of the partition is at position $(0,0)$, the positions in the first row go up to $(0,k-1)$, and the next row starts at $(1,0)$ and proceeds from there.

A reader unfamiliar with abacus techniques may wish to prove the first and second facts of Lemma \ref{CoreFacts} as a short exercise using these definitions.

\vspace{0.1in}

A \emph{numerical set} $S$ is a subset of $\mathbb{N} \bigcup \{0\}$ which contains 0 and has finite complement.  If $S$ is additively closed, i.e. $a+b \in S \, \forall \, a,b \in S$, then $S$ is a \emph{numerical monoid} or \emph{numerical semigroup}; both terms are used in the literature, and since we assume $0 \in S$ we use the former.

The number of elements in the complement of a numerical monoid is its \emph{genus} $g$, and the largest element of the complement is its \emph{Frobenius number} $f$.  Any numerical set contains an \emph{atom monoid} $A(S)$, the largest subset of $S$ closed under translation by $S$: $A(S) = \{n \in \mathbb{N} \bigcup \{0\}\vert n+s \in S\ \forall s \in S \}$.  We know $A(S) \subseteq S$ because $0 \in S$.  Clearly if $S$ is a numerical monoid then $S = A(S)$.

The Frobenius number of the atom monoid is the same as the largest element of the complement of $S$, since $f \not\in S \Rightarrow f+0=f \not\in A(S)$ and $f+c \in S$ for all $c>0$. (For more results on numerical sets, see Marzuola and Miller, \cite{MarzMiller}).  Thus, in the sense of largest possible complement, the minimum possible atom monoid for a numerical set with Frobenius number $f$ is $\{0,f+1,f+2,f+3,\dots\}$.

\vspace{0.1in}

\noindent \textbf{Example:} The set $S = \{0,3,6,7,9, 12\} \bigcup \{n \in \mathbb{N} : n > 12\}$ is a numerical set.  Its atom monoid is $A(S) = \{0, 6, 9, 12\} \bigcup \{n \in \mathbb{N} : n > 12\}$.  The genus of $A(S)$ is 9; its Frobenius number is 11.

\section{Partitions and Numerical Monoids}

The bijection $\phi:S \rightarrow \lambda$ which proves Theorem \ref{HooksetMonoids} is as follows.

\vspace{0.1in}

\noindent \textbf{The bijection:} Build the profile of a partition $\lambda$ by numbering the steps east and north from the bottom left corner starting with 0 for the first east step.  Given a numerical set $S$ with Frobenius number $f$, for $0 < i \leq f$, let step $i$ be east if $i \in S$ and north if $i \not\in S$.  We may think of concluding the partition profile with an infinite horizontal line, if desired, but parts are counted only between the nontrivial steps.  The map is plainly bijective; the reverse process from $\lambda$ to $S$ is to read off included and excluded elements with each step.

\vspace{0.1in}

\noindent \textbf{Example.} Consider the partition $\lambda = (5,4,4,2)$:

\begin{center}\begin{tabular}{c} $\Box \Box \Box \Box \Box$\vspace{-0.05in} \\ $\Box \Box \Box \Box$\phantom{$\Box$}\vspace{-0.05in}\\$\Box \Box \Box \Box$\phantom{$\Box$}\vspace{-0.05in}\\$\Box \Box$\phantom{$\Box \Box \Box$}\vspace{-0.05in}\end{tabular}\end{center}

Assign 0 to the first east step.  Then the north steps occur at steps 2, 5, 6, and 8.  The numerical set associated to $\lambda$ by $\phi$ is thus given by $\{0,1,\, \, 3,4, \, \, 7, \, \, 9, 10, \dots \}$.

Corollary \ref{EmptyAtom}, as reworded from Olsson's Theorem 2.2, follows from our Theorem \ref{HooksetMonoids}, and Theorem \ref{HooksetMonoids} is a result of the following property of this map:

\begin{theorem}\label{PtnAtomMonoid} The map $\phi$ satisfies $Hk(\lambda) = \overline{A(S)}$, that is, it maps $k$-cores to numerical sets with $k$ in the atom monoid and vice versa.
\end{theorem}

\noindent \textbf{Proof:} Hook lengths in $\lambda$ are precisely the differences between the numerical identifier of north steps, and east steps of lower numerical identifier.

Suppose that $k$ is a hook length in $Hk(\lambda)$, i.e., $\lambda$ is not a $k$-core.  Then there is some east step at a place $i$, hence an $i \in S$, with a north step (bead) at place $i+k$, so $i+k \not\in S$.  Then by the definition of $A(S)$, $k \not\in A(S)$.

Vice versa, suppose $k$ is not a hook length in $Hk(\lambda)$.  Then for any $i \in S$, $i + k \in S$, and so $k \in A(S)$. $\Box$

\vspace{0.1in}

The Introduction mentioned that a similar map can be found in one earlier paper.  The map in \cite{BAdMier}, by Maria Bras-Amor\'{o}s and Anna de Mier, makes the reverse assignment of vertical and horizontal steps.  It starts the path at the origin and, the genus $g$ having been specified, bounds the path within a square box of size $g$ in the first quadrant.  That map has the virtue that a statistic of interest in semigroup studies, the weight of the semigroup, is equal to the weight of the partition bounded within a square above the path.  To give it the hook property of Theorem \ref{PtnAtomMonoid}, one may adjust the shape of the upper boundary by moving the left side of their bounding square back one unit and counting the hook sizes of the partition then bounded above the path.  The bounding square otherwise becomes irrelevant for hook-counting, which was not a topic of interest in that paper (which was primarily concerned with the symmetry of the path, and the Dyck partition bounded below the path).

\subsection{Implications for partitions}

The fact that hooksets are the complements of additive monoids can extend previous work on the hooksets of partitions.  Anderson \cite{Anderson} showed

\begin{theorem} If $s$ and $t$ are coprime, the number of partitions which are simultaneously $s$-core and $t$-core is finite. \end{theorem}

Olsson and Stanton \cite{STCoresPrime} explicitly found the largest such partition, showing that the largest simultaneous $s$-core and $t$-core partitions $(s^2-1)(t^2-1)/24$.

Let $S = \{ s_1, \dots , s_k \}$ be any set of positive integers.  Call \emph{$S$-core} the set of partitions that are simultaneously $s_i$-core for all $s_i \in S$.  The correspondence above immediately yields Theorem \ref{SCores}: \vspace{0.1in}

\noindent \textbf{Proof:} Clearly if $gcd(S) = c > 1$, any $c$-core is simultaneously $s_i$-core for all $s_i$, and these sets are known to be infinite.

Assume $gcd(S) = 1$.  Being the complement of a numerical monoid, if two hooklengths $s$ and $t$ are missing from the hookset of a partition, then so is hooklength $s+t$.  (Direct proof: if $s$ is missing from the hookset, any spacer in the abacus must be followed, $s$ places later, by another spacer; that spacer must be followed $t$ places later by another spacer.)  Since a set $S$ for which $gcd(S)=1$ has finite Frobenius number $f$, the set of missing hooks includes all integers greater than $f$, and so the outer hook of any $S$-core partition is at most $f$.  There are $2^{f-1}$ such partitions. \hfill $\Box$

\vspace{0.1in}

This also gives the crude bound $\left( \frac{f}{2} \right)^2$, i.e. a square of outer hook $f$, for the largest $n$ partitioned by an $S$-core.  Necessarily, the actual $n$ will be significantly less.  Arguments on the profile shape can reduce this considerably for any particular values in $S$, and such partitions will generally be quite "concave."  For example, square partitions have all hook sizes, but avoiding hooks of size $k$ requires that any segment of the profile of $k$ units length cannot be straight.

For a set of two coprime elements $S=\{s,t\}$, if $s$ and $t$ are coprime then Sylvester showed that the Frobenius number is $f = st-s-t$.  This gives us $n < \frac{(st-s-t)^2}{4}$, which for large $s$ and $t$ is about 6 times the size of Olsson and Stanton's $\frac{(s^2-1)(t^2-1)}{24}$.

\vspace{0.1in}

\noindent \textbf{Remark:} It should also be noted that in the proof above, finiteness of the set of partitions when $gcd(S)=1$ is the main interest for the connection with numerical monoids, and infiniteness for $gcd(S)=g>1$, $g \not\in S$, is simply proved by referring to $g$-cores.  However, Aukerman et al. prove in \cite{STCoresNotPrime} that the number of simultaneous $s/t$-cores that are not $g$-cores is infinite as well, and this we can also do, using an explicit abacus construction outlined in \cite{Nath}.  Suppose $gcd(S)=g>1$, with $g \not\in \{S\}$.  Say $max(S)=s$.  Build an $s$-abacus by filling columns $s-1$, $(s-1)-g$, $(s-1)-2g$, $\dots$ to height $m \geq 0$ with beads, and remove from the top row any selection of beads that satisfies the following conditions:

\begin{itemize}
\item at least one bead remains with a bead removed (leaving a spacer) $g$ places below;
\item no bead is removed $s_i$ places (that is, $(s_i/g) \cdot g$ places) below a bead that remains, for any $s_i$.
\end{itemize}

Then the remaining abacus is that of an $S$-core, since it was originally a $g$-core and no hook of size $(s_i/g)\cdot g$ is formed by removals.  But the first condition gives us a $g$-hook, so the partition is not a $g$-core.  For any given permissible arrangement of beads on the top row, the height can be extended indefinitely, giving us an infinity of such partitions.

\vspace{0.1in}

\noindent \textbf{Example:} Let $S=\{8,12\}$.  On a 12-abacus, fill columns 11, 7, and 3 to some equal height.  Remove the bead in column 7, creating a 4-hook.  (We cannot remove 3 without removing 11, since that would create an 8-hook; we could remove both 3 and 11.)  The abacus will be

\begin{center}\begin{tabular}{c}
$\circ \circ \circ \bullet \circ \circ \circ \bullet \circ \circ \circ \bullet$ \vspace{-0.07in} \\
\hline $\circ \circ \circ \bullet \circ \circ \circ \bullet \circ \circ \circ \bullet$ \vspace{-0.085in} \\
\hline \dots \\
$\circ \circ \circ \bullet \circ \circ \circ \bullet \circ \circ \circ \bullet$ \vspace{-0.07in} \\
\hline $\circ \circ \circ \bullet \circ \circ \circ \circ \circ \circ \circ \bullet$ \vspace{-0.075in} \\ \hline
\end{tabular}\end{center}

With a little more verbiage, one can give the procedure more flexibility by working in each of the residue classes of columns mod $g$, and allowing removed beads to dip down into rows past the first, though the second condition restricts how far down removals can go in each residue group.  Different residue groups mod $g$ may be raised to independent heights, although the group for residue 0 may not be raised indefinitely since the 0 column must be empty.

In fact, this set of small $S$-cores being borne aloft on growing $g$-cores is the combinatorial meaning of the Aukerman/Kane/Sze product between the infinite generating function of $g$-cores, and the polynomial generating function of the $\frac{s}{g}/\frac{t}{g}$-cores.  The non-$g$-cores are simply those which begin with a non-$g$-core $\frac{s}{g}/\frac{t}{g}$-core; one such as the example with which we began will always exist.

\subsection{Implications for numerical sets}

This map might also be fruitful in considering some of the open questions in numerical monoid theory.

Theorem \ref{PtnAtomMonoid} was originally motivated by an interest in partitions with hooksets that are an initial segment of the integers, i.e. $Hk(\lambda) = \{1,2,3,\dots,g\}$.  One might call these \emph{non-cores}, since they are $k$-core for the fewest possible $k$: any partition is a $k$-core for $k$ larger than its outer hook, and these are not $k$-core for any smaller $k$.  These were counted by computation, and the sequence counting non-cores with largest outer hook $g$ was located in Neil Sloane's Online Encyclopedia of Integer Sequences \cite{SloaneSeq}, where the sequence appeared as number A158291, titled "The number of numerical sets $S$ with atom monoid $A(S)$ equal to $\{0,n+1,n+2,n+3,n+4,...\}$."  The map above was then quickly produced.

Since a numerical monoid is its own atom monoid, the partitions which map to numerical monoids form representatives of the sets of partitions with particular hooksets.  It is known (see Chung and Herman, \cite{ChungCollab}) that neither the hookset nor even the hook multiset specifies the partition; Chung and Herman calculate a few of the values of the number of partitions with equal hooksets, and ask for a general formula.

An open conjecture regarding numerical monoids (see \cite{BAFibonacci}) is that $n_g$, the number of numerical monoids with genus $g$, grows in roughly Fibonacci fashion, i.e. $n_g \approx n_{g-1} + n_{g-2}$.  Chung and Herman's formula, for the number of partitions with outer hook $g$ that are not specified by their hooksets, would be just the difference between $2^{g-1}$ -- which is the number of partitions with outer hook $g$ -- and the number $n_g$.

An alternative proof strategy for this conjecture involving partitions might be to use a growth process which constructs new hooksets from old by adding new sizes of hook in a fashion that typically adds 1 or 2 at a time, in a mostly reversible way.  There is a variety of such growth processes in the partition literature that might be considered.

\section{Partitions non-core for all but one $k$}

As mentioned in Lemma \ref{CoreFacts}, the only 2-cores are the staircase partitions, which have hooksets $Hk(\lambda) = \{1,3,5,\dots,2k-1\}$.  Therefore all 2-cores with largest hook $2k-1$ possess as hooks all nonmultiples of 2 from 1 to their maximum.

If $k > 2$, it is possible to have $k$-cores with hooksets that are contained within the nonmultiples of $k$, but miss some values between the 1 and the largest nonmultiple of $k$ appearing.  For instance, $Hk((3,1,1)) = \{1,2,5\}$, which is a 3-core in which the nonmultiple 5 appears as a hooklength, but 4 does not.  Likewise, a 2-core is a 4-core that is also missing the hooks which are congruent to 2 mod 4.

Here we ask for behavior similar to having a hookset which is an initial segment of the integers, with one exception: what characterizes $k$-cores with complete hooksets, i.e. $Hk(\lambda)$ is an initial segment of the nonmultiples of $k$?  What is the generating function for these partitions?

Several of the argument techniques may be illustrated if we begin with the simplest next case, that of 3-cores.

\begin{theorem} The generating function of the 3-cores with complete hooksets is $$ \sum cc_3 (n) q^n = 1+q + 2 \left( \sum_{k \geq 2} q^{k^2}+q^{k^2-k} \right)  \, \text{ .}$$
\end{theorem}

\noindent \textbf{Proof:} The abaci of 3-cores consist of uninterrupted columns of beads on 3 rods, which we will call places 0, 1, and 2.  (So a column is the \emph{occupied} portion of a place.)  Abaci representing 3-cores with all possible hooks are as follows:

\begin{center}
\begin{tabular}{c}
$\circ \bullet \circ$ \vspace{-0.07in} \\
\hline $\circ \bullet \circ$ \vspace{-0.085in} \\
\hline \dots \\
$\circ \bullet \circ$ \vspace{-0.075in} \\ \hline
\end{tabular} \quad
 \begin{tabular}{c}
$\circ \circ \bullet$ \vspace{-0.07in} \\
\hline $\circ \circ \bullet$ \vspace{-0.085in} \\
\hline \dots \\
$\circ \circ \bullet$ \vspace{-0.075in} \\ \hline
\end{tabular} \quad
\begin{tabular}{c}
$\circ \bullet \bullet$ \vspace{-0.07in} \\
\hline $\circ \bullet \bullet$ \vspace{-0.085in} \\
\hline \dots \\
$\circ \bullet \bullet$ \vspace{-0.075in} \\ \hline
\end{tabular} \quad
\begin{tabular}{c}
$\circ \bullet \bullet$ \vspace{-0.07in} \\
\hline $\circ \bullet \bullet$ \vspace{-0.085in} \\
\hline \dots \\
$\circ \bullet \circ$ \vspace{-0.075in} \\ \hline
\end{tabular}\end{center}

The first abacus illustrates a partition with one column of beads down the middle.  If the last row with a bead occupying it is the $k$-th row, this represents a partition $(1,3,5,7,...,2k-1)$.  Its outer hook will be of length $3k-2$.

Hooks of all lengths not multiples of 3 can be read from the bottommost bead to the spacers above it in the columns to its left and right, from 1 up to the outer hook.  For instance, the distance from the bottommost bead to the spacer immediately to its left is a hook of length 1; from the bottommost bead to the spacer immediately above and to its right we have a hook of length 2.  The spacer immediately above and left of the bottommost bead gives a hook of length 4, and we continue up the abacus.  Thus, it holds all possible hooks that are not multiples of 3.

The second abacus represents $(2,4,6,8,...,2k)$.  Its hooks can be read from the bottommost bead to the columns to its left.

The third abacus represents a 3-core with both places 1 and 2 occupied by columns of equal height.  This is the partition $(1,1,2,2,3,3,...,k,k)$.  The fourth abacus represents a 3-core with the column in place 2 of length 1 bead shorter than that in place 1, so that the last bead is in place 1.  This is the partition $(1,1,2,2,3,3,...,k)$.

In the latter two abaci, not all the possible hooks can be read from the last bead.  For example, in the third abacus, the hooks which would have been of length 1 mod 3 if read from the last bead in place 2 would require spacers in place 1, but these are all occupied by beads.  However, from those beads, the hooks of length 1 mod 3 can be read from place 1 to the always-empty place 0.

Thus, imagine place 2 to have been occupied first.  Then causing place 1 to be occupied destroyed some hooks as read from the last bead, but replaced those destroyed.  When we add such a column, the hooks created must be of lengths that replace all of those destroyed, especially the longest.  This will be the mental model we use in larger abaci.

On the other hand, in the remaining possible types of 3-core abacus, namely 

\begin{center}
\begin{tabular}{c}
$\circ \bullet \bullet$ \vspace{-0.07in} \\
\hline $\circ \bullet \bullet$ \vspace{-0.085in} \\
\hline \dots \\
$\circ \circ \bullet$ \vspace{-0.075in} \\ \hline
\end{tabular} \quad
\begin{tabular}{c}
$\circ \bullet \bullet$ \vspace{-0.07in} \\
\hline $\circ \bullet \circ$ \vspace{-0.085in} \\
\hline \dots \\
$\circ \bullet \circ$ \vspace{-0.075in} \\ \hline
\end{tabular}
\end{center}

\noindent with the column in place 1 shorter than that in place 2, or that in place 2 at least 2 beads shorter than that in place 1, we miss at least one hook.

For the left abacus, the longest column is place 2, so think of this as having been placed first.  The largest hook of length 1 mod 3 from the last bead is destroyed (there is a bead at the bottom of place 1 now), but not replaced by the column added in place the column in place 1 (it would have to be of equal height).  For the right abacus, the largest hook of length 2 mod 3, as read from the last bead in place 1, is destroyed by adding a column in place 2, and not replaced.  Thus, the first four abaci above are all possible 3-cores with complete hooksets, and their generating function by weight is as given.  $\Box$

This classification of the behaviors of columns in the abaci of $k$-cores with complete hooksets will persist.  Columns will be one of four types: the last longest column; necessarily empty; necessarily of height equal to the longest (or 1 less, if in a higher place); or, once $k \geq 5$, fillable to any height between 0 and the maximum without losing any hooks.

\subsection{Construction Algorithm}

Given $k$, an algorithm to produce the abacus of a $k$-core with a hookset which is the initial segment of the whole numbers not multiples of $k$ is as follows:

\begin{itemize}
\item 1) Choose a place $p_1$ and a height $m$ to be the highest (rightmost) place among those filled to the maximum height $m$.
\item 2) Choose a set of places $\{p_i\}$ to be filled to the maximum possible height: $m$ if $p_i \leq p_1$, and $m-1$ if $p_i > p_1$.  The condition constraining this choice is as follows: construct the set of differences of $\{p_i\}$ and its complement, but only include negative differences $p_i - p_j$, $p_j > p_i$, if $p_i \leq p_1$, i.e. the set $\{p_a - p_b \vert p_a \in \{p_i\}, p_b \in \overline{\{p_i\}}, p_b > p_a \, \text{only if} \, p_i \leq p_1 \}$. This difference set must include at least one representative of all possible nonzero residues modulo $k$.
\item 3) Construct a new set $\{p_r\}$ as follows.  Consider the remaining places $p_a \in \overline{\{p_i\}}$.  The spacer $p_a$ is \emph{redundant} if the distance $p_1 - p_a$ is repeated by a bead $p_j$ and a spacer $p_b$ with $p_j < p_1$, or with $p_a > p_1$, $p_j > p_1$, and $p_j > p_b$.  If desired, add $p_a$ to $\{p_r\}$ and repeat Step 3 for the other spacers in $\overline{\{p_r\} \bigcup \{p_i\}}$, considering only distances to spacers not in $\{p_r\}$.  Elements chosen for $\{p_r\}$ may be filled to any height up to 1 less than the maximum: $m-1$ if $p_a < p_1$, and $m-2$ if $p_a > p_1$.
\end{itemize}

\noindent \textbf{Example:} Let $k=5$, $m>1$.  We wish to build 5-cores with all possible hooks that are nonmultiples of 5.  We can choose columns from among places 1 through 4.

In step 1, we are free to choose any place.  Suppose we choose place 2.  Places 1, 3, and 4 are currently open; the respective residues represented by the distances from the bead in place 2 to the spacers in places 0, 1, 3, and 4 are 2, 1, -1, and -2 respectively, which are all of the nonzero residues modulo 5.

In step 2, we can choose a set of places that leaves the remaining residues as a full list, under the conditions prescribed.  For example, we could choose to fill place 1 to height $m$ or place 3 or place 4 to height $m-1$.  The last choice makes $\{p_i\} = \{2,4\}$.  Then spacers would occupy places 0, 1, and 3, and the differences would be 2, 1, and -1 (from the beads in place 2) and 4, 3, and 1 (from the beads in place 4).

We could not choose $\{p_i\} = \{2,3\}$, since with beads in places 2 and 3, and spacers in places 0, 1, and 4, the distances to spacers would be 2, 1 and -2 (from the bead in place 2), and 2 and 3 (from the bead in place 3, only counting negatives).  For instance, if we fill place 2 to height 2 and place 3 to height 1, there is no hook of length 4.

Let us say we choose $\{p_i\} = \{2,4\}$.  What spacers are redundant?  Place 1 is not, since 4-3 is not a permissible replacement for 2-1.  (The hook of length $5m+1$ would not be replaced from place 4, which is of height $m-1$.)  Place 3 is redundant, since 2-3 is replaced by 4-0, and this is permissible.  (The hook of length $5m-1$ from the top bead is replaced by that hook from the top of place 4.)

A brief remark on some variants that might be illustrative.  If we had not added place 4 to $\{p_j\}$, place 3 would not have been a redundant spacer.  If on the other hand $k$ were much larger, and place $k-1$ were our $p_1$, then place $k-4$, place 1 in this example, \emph{would} be redundant, since the hook lost would be replaced from the lower column, place $k-3$.  It is also possible, in a larger case, to have one choice of spacer make another no longer redundant.

In this case, we have specified that place 2 is to be built to height $m$, place 4 is to built to height $m-1$, place 1 will remain empty, and place 3 will have some content from 1 up to $m-2$.  These particular 5-cores will have as hooks all nonmultiples of 5 from 1 up to $5m-3$.  The abacus will look similar to this:

\begin{center}
\begin{tabular}{c}
$\circ \circ \bullet \bullet \bullet$ \vspace{-0.07in} \\
\hline $\circ \circ \bullet \bullet \bullet$ \vspace{-0.07in} \\
\hline \dots \\
$\circ \circ \bullet \bullet \bullet$ \vspace{-0.07in} \\
\hline $\circ \circ \bullet \circ \bullet$ \vspace{-0.07in} \\
\hline $\circ \circ \bullet \circ \circ$ \vspace{-0.075in} \\ \hline
\end{tabular}
\end{center}

If we fill place 2 to height $m$, and place 3 to height $0 \leq c \leq m-2$, the generating function for the 5-cores so specified is

$$ \sum_{m=1}^{\infty} \sum_{c=0}^{m-2} q^{3(c^2 + c)} q^{(2c)(2(m-c)-1)} q^{(m-c)^2 + (m-c) + \frac{3}{2}\left( (m-c-1)^2+(m-c-1) \right)} \, \text{.}$$

Running over the various permissible top-level profiles and summing, we would obtain the generating function for 5-cores with complete hooksets.

\section{Arbitrary Hooksets}

Less can be said about less-structured conditions on hooksets.  Given an arbitrary set $S$ of hooks which we require the partition $\lambda$ to contain, a sufficient condition (by no means necessary) for $Hk(\lambda)$ to contain $S$ can be established by the following argument on the abacus:

Duplicate the unwrapped abacus of $\lambda$ and lay the two copies side by side, with the lower copy offset to the left by the distance $s_i$.  If any spacer in the top copy has a bead below it, there exists somewhere in the partition a hook of length $s_i$.

\begin{center}\begin{tabular}{c}
\includegraphics[scale=0.8]{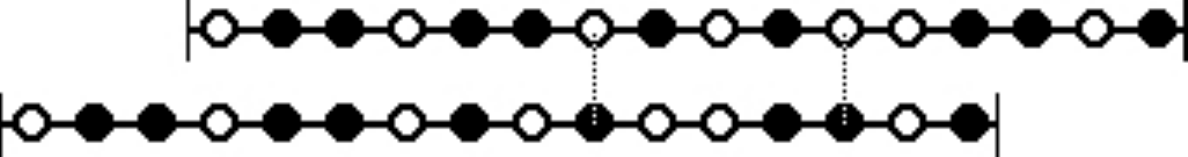} \\
Positions 6 and 10 (counted from 0) indicate the\\
 presence in this partition of hooks of length 3.
\end{tabular}\end{center}

For an offset of distance $s_i$, the overlap is of length equal to that of the profile of the partition, less $s_i$.  The overlapping part of the top copy represents the profile of the partition up to a distance $s_i$ steps from its uppermost corner, and likewise the bottom part from the upper corner down to $s_i$ steps from the lowermost left.

If the profile of the partition lies strictly above the diagonal at this point on the upper copy, it means that spacers populate more than half the places in the overlapping section from the start of the profile, and likewise the bottom portion being above the diagonal means beads populate more than half the places in the overlapping section from the upper end.  The Pigeonhole Principle then tells us that somewhere in the overlap must be a match of a spacer to a bead.  This proves Theorem \ref{DyckHooks}.

It follows as Corollary \ref{StrictCorollary} that strict Dyck paths -- meaning that the profiles lie everywhere strictly above the diagonal except at the corners, equal distances from the origin -- are a subset of the partitions discussed earlier which possess a hookset that is an initial segment of the integers.  (This is not precisely a new result; using a different methodology it can be proved with little difficulty as a follow-up to problem 7.107 of Richard Stanley's Enumerative Combinatorics, Vol. 2.)


\begin{thebibliography}{99}
\bibitem{Anderson} Anderson, J.  Partitions which are simultaneously $t_1$-core and $t_2$-core.  Discrete Math. 248 (2002) 237-243
\bibitem{STCoresNotPrime} Aukerman, D., Kane, B., and Sze, L.  On simultaneous $s$-cores/$t$-cores.  Discrete Mathematics 309 (2009) 2712-2720
\bibitem{BAFibonacci} Bras-Amor\'{o}s, M. Fibonacci-like behavior of the number of numerical semigroups of a given genus.  Semigroup Forum (2008) 76: 379Ð384.  DOI 10.1007/s00233-007-9014-8
\bibitem{BAdMier} Bras-Amor\'{o}s, M., and de Mier, Anna.  Representation of numerical semigroups by Dyck paths.  Semigroup Forum, Springer, vol. 75, n. 3, pp. 676-681, December 2007. ISSN: 0037-1912.  arXiv: http://arxiv.org/abs/math/0612634
\bibitem{RepSource} Brauer, R., and Robinson, G. de B.  On a conjecture by Nakayama.  Trans. Roy. Soc. Canada Sect. III. (3) 41 (1947), 11-25
\bibitem{ChungCollab} Chung, F. R. K., and Herman, J. E.  Some results on hook lengths.  Discrete Mathematics 20 (1977) 33-40
\bibitem{CrankSource} Garvan, F.G., Kim, D., and Stanton, D.  Cranks and $t$-cores.  Inventiones Mathematicae 101 (1990), 1-17
\bibitem{MarzMiller} Marzuola, J., and Miller, A.  Counting numerical sets with no small atoms.  http://arxiv.org/abs/0805.3493
\bibitem{Nath} Nath, R. On the s-core of a t-core partition, Integers: Electronic Journal Of
Combinatorial Number Theory 8 (2008)
\bibitem{STCoresPrime} Olsson, J. B., and Stanton, D.  Block inclusions and cores of partitions.  Aequat. Math. 74 (2007), 90-110
\bibitem{AddArt} Olsson, J.B.  Core partitions and block coverings.  Proc. Amer. Math. Soc. 137 (2009), 2943-2951
\bibitem{SloaneSeq} Online Encyclopedia of Integer Sequences, Sequence A158291.  http://www.research.att.com/~njas/sequences/A158291
\end{thebibliography}
\end{document}